\documentclass[11pt]{article}%
\usepackage{amsmath}%
\usepackage{amsfonts}%
\usepackage{amssymb}%
\usepackage{graphicx}
\usepackage{amsthm}
\newtheorem{theorem}{Theorem}

\newtheorem{proposition}[theorem]{Proposition}
\theoremstyle{definition}
\newtheorem{definition}[theorem]{Definition}

\newtheorem{remark}[theorem]{Remark}
\newtheorem{example}[theorem]{Example}

\begin{document}

\title{On the Tensor Permutation Matrices}

\author{Christian Rakotonirina}
\maketitle

\noindent D\'{e}partement du G\'{e}nie Civil
\\ Institut Sup\'{e}rieur de Technologie d'Antananarivo, IST-T, BP 8122\\
\noindent D\'{e}partement de Physique
\\ Laboratoire de Rh\'{e}ologie des Suspensions, LRS, Universit\'{e} d'Antananarivo
\\ Madagascar.\\
\noindent e-mail : rakotopierre@refer.mg

\begin{abstract}
\noindent We show that two TPM's  permute tensor product of rectangle matrices. An example, in the particular case of tensor commutation matrices (TCM's), for studying a linear matrix
equation is given.
\emph{\textbf{MCS 2010}}:  15A69
 \end{abstract}

\section*{Introduction}

When we were working on Raoelina Andriambololona idea on the using
tensor product of matrices in the Dirac equation \cite{rak03}, \cite{wan01}, we
met the unitary
matrix\\
\begin{center}
$\textbf{U}_{2\otimes2}=
\begin{bmatrix}
  1 & 0 & 0 & 0 \\
  0 & 0 & 1 & 0 \\
  0 & 1 & 0 & 0 \\
  0 & 0 & 0 & 1 \\
\end{bmatrix}$
\end{center}
which has the following properties: for any unicolumns and two rows
matrices\\
$\textbf{a}=
\begin{bmatrix}
  \alpha^1 \\
  \alpha^2 \\
\end{bmatrix}\in\mathbb{C}^{2\times1}$,
$\textbf{b}=
\begin{bmatrix}
  \beta^1 \\
  \beta^2 \\
\end{bmatrix}\in\mathbb{C}^{2\times1}$
\begin{equation*}
\textbf{U}_{2\otimes2}\cdot(\textbf{a}\otimes\textbf{b})=\textbf{b}\otimes\textbf{a}
\end{equation*}
and for any two $2\times2$-matrices, $\textbf{A}$,
$\textbf{B} \in\mathbb{C}^{2\times2}$
\begin{equation*}
\textbf{U}_{2\otimes2}\cdot(\textbf{A}\otimes\textbf{B}) = (\textbf{B}\otimes\textbf{A})\cdot\textbf{U}_{2\otimes2}
\end{equation*}
This matrix is frequently found in quantum information theory
\cite{fuj01}, \cite{fad95}, \cite{fra02}.

\noindent We call this matrix a tensor
commutation matrix (TCM) $2\otimes2$. The TCM
$3\otimes3$ has been written by Kazuyuki Fujii \cite{fuj01} under the
following form\\

\begin{equation*}
\textbf{U}_{3\otimes3}=
\begin{bmatrix}
 
\begin{bmatrix}
  1 & 0 & 0 \\
  0 & 0 & 0 \\
  0 & 0 & 0 \\
\end{bmatrix}
  
\begin{bmatrix}
  0 & 0 & 0 \\
  1 & 0 & 0 \\
  0 & 0 & 0 \\
\end{bmatrix}
  
\begin{bmatrix}
  0 & 0 & 0 \\
  0 & 0 & 0 \\
  1 & 0 & 0 \\
\end{bmatrix}
 \\
 
\begin{bmatrix}
  0 & 1 & 0 \\
  0 & 0 & 0 \\
  0 & 0 & 0 \\
\end{bmatrix}
  
\begin{bmatrix}
  0 & 0 & 0 \\
  0 & 1 & 0 \\
  0 & 0 & 0 \\
\end{bmatrix}
  
\begin{bmatrix}
  0 & 0 & 0 \\
  0 & 0 & 0 \\
  0 & 1 & 0 \\
\end{bmatrix}
 \\
  
\begin{bmatrix}
  0 & 0 & 1 \\
  0 & 0 & 0 \\
  0 & 0 & 0 \\
\end{bmatrix}
  
\begin{bmatrix}
  0 & 0 & 0 \\
  0 & 0 & 1 \\
  0 & 0& 0 \\
\end{bmatrix}
  
\begin{bmatrix}
  0 & 0 & 0\\
  0 & 0 & 0 \\
  0 & 0 & 1 \\
\end{bmatrix}
 
\end{bmatrix}
\end{equation*}
in order to obtain a conjecture of the form of a TCM $n\otimes n$, for any $n$ $\in$ $\mathbb{N}^\star$. He calls these matrices "`swap operator"'.\indent

  $\textbf{U}_{n\otimes
p}$, the TCM $n\otimes p$, $n$,
$p$\;$\in$\;$\mathbb{N}^\star$, commutes the tensor product of $n\times n$- matrix by $p\times p$-matrix. In this paper we will show that two $\sigma$-TPM's $\textbf{U}_{\sigma}$, $\textbf{V}_{\sigma}$   permute tensor product of rectangle matrices, that is, $\textbf{U}_{\sigma}\cdot\left(\textbf{A}_1\otimes \textbf{A}_2\otimes
\ldots \otimes\textbf{A}_k\right)\cdot\textbf{V}_{\sigma}^T=[\textbf{A}_{\sigma(1)}\otimes
\textbf{A}_{\sigma(2)}\otimes \ldots\otimes
\textbf{A}_{\sigma(k)}$, where $\sigma$ is a permutation of the set
$\left\{1, 2,\ldots , k\right\}$.\\
 $\textbf{U}_{\sigma}=\textbf{V}_{\sigma}$, if $\textbf{A}_1$, $\textbf{A}_2$, $\ldots$, $\textbf{A}_k$ are square matrices (Cf. for example $\cite{rak05}$).\\
We will show this property, according to the Raoelina Andriambololona approach in linear
and multilinear algebra \cite{rao86}: in establishing at first, the propositions on linear operators in intrinsic way, that is
independently of the bases, and then we demonstrate the
analogous propositions for the matrices.

 \section*{Tensor Product of Matrices}
\begin{definition}
Consider $\textbf{A}=(A^i_j)\in\mathbb{C}^{m\times n}$, $\textbf{B}=(B^i_j)\in\mathbb{C}^{p\times r}$. The matrix defined by
\begin{equation*}
\textbf{A}\otimes\textbf{B}=
\begin{bmatrix}
  A^1_1\textbf{B} & \ldots & A^1_j\textbf{B} & \ldots & A^1_n\textbf{B} \\
  \vdots &  & \vdots &  & \vdots \\
  A^i_1\textbf{B} & \ldots & A^i_j\textbf{B} & \ldots & A^i_n\textbf{B} \\
  \vdots &  & \vdots &  & \vdots \\
  A^m_1\textbf{B} & \ldots & A^m_j\textbf{B} & \ldots & A^m_n\textbf{B} 
\end{bmatrix}
\end{equation*}
obtained after the multiplications by
scalar, $A^i_j\textbf{B}$, is called the tensor product of the matrix
$\textbf{A}$ by the matrix $\textbf{B}$.
\begin{equation*}
\textbf{A}\otimes \textbf{B}\in \mathbb{C}^{mp\times
nr}
\end{equation*}
\end{definition}
\begin{proposition}
Tensor product of matrices is associative.
\end{proposition}
\begin{proposition}
Consider the linear operators $A$\;$\in$\;$\mathcal{L}( \mathcal{E},
\mathcal{F} )$ , $B$\;$\in$\;$\mathcal{L}( \mathcal{G}, \mathcal{H}
)$. $\textbf{A}$ is the matrix of A with respect to the couple of bases
$\left(%
  (\overline{e_i})_{1\;\leq i\;\leq n}, (\overline{f_j})_{1\;\leq
j\;\leq m} \\
\right)$, $\textbf{B}$  the one of B in
$\left(%
(\overline{g_k})_{1\;\leq
k\;\leq r}, (\overline{h_l})_{1\;\leq l\;\leq p}\\
\right)$. Then, $\textbf{A}\otimes \textbf{B}$ is the matrix of
$A\otimes B$ with respect to the couple of bases ($\mathcal{B}$,
$\mathcal{B}_1$), where\\
$\mathcal{B}$\;=\;$(\overline{e_1} \otimes \overline{g_1},
\overline{e_1} \otimes \overline{g_2}, \ldots ,\overline{e_1}
\otimes \overline{g_r}, \overline{e_2} \otimes \overline{g_1},
\overline{e_2} \otimes \overline{g_2}, \ldots , \overline{e_2}
\otimes \overline{g_r}, \ldots , \overline{e_n} \otimes
\overline{g_1}, \overline{e_n} \otimes
\overline{g_2}, \ldots , \overline{e_n} \otimes \overline{g_r})$\\
\end{proposition}
\noindent Notation: we denote the set $\mathcal{B}$ and $\mathcal{B}_1$ by\\
\begin{equation*}
 \mathcal{B}=\left(%
 \overline{e_i} \otimes \overline{g_k}
\right)_{1\;\leq i\;\leq n, 
        1\;\leq k\;\leq r}
        =\left(%
  (\overline{e_i})_{1\;\leq i\;\leq n} 
\right) \otimes \left(%
(\overline{g_k})_{1\;\leq
k\;\leq r}
\right)
\end{equation*}

\begin{equation*}
 \mathcal{B}_1=\left(%
 \overline{f_j}\otimes \overline{h_l}
\right)_{1\;\leq j\;\leq m, 
        1\;\leq l\;\leq p}\\
        =\left(%
  (\overline{f_j})_{1\;\leq j\;\leq m} 
\right) \otimes \left(%
(\overline{h_l})_{1\;\leq
l\;\leq p}
\right)
\end{equation*}

\section*{Tensor permutation operators}
\begin{definition}
Consider the  $\mathbb{C}$- vector spaces $\mathcal{E}_1$,
$\mathcal{E}_2$, \ldots, $\mathcal{E}_k$ and a permutation $\sigma$ of
$\left\{1, 2,\ldots , k\right\}$. The linear operator $U_{\sigma}$ 
   from $\mathcal{E}_1$$\otimes$
$\mathcal{E}_2$$\otimes$\ldots$\otimes$ $\mathcal{E}_k$ to
$\mathcal{E}_{\sigma(1)}$$\otimes$
$\mathcal{E}_{\sigma(2)}$$\otimes$\ldots$\otimes$
$\mathcal{E}_{\sigma(k)}$, $U_{\sigma}$$\in$
$\mathcal{L}$($\mathcal{E}_1$$\otimes$
$\mathcal{E}_2$$\otimes$\ldots$\otimes$ $\mathcal{E}_k$,
$\mathcal{E}_{\sigma(1)}$$\otimes$
$\mathcal{E}_{\sigma(2)}$$\otimes$\ldots$\otimes$
$\mathcal{E}_{\sigma(k)}$ ), defined by
\begin{center}
$U_{\sigma}$($\overline{x_1}$\;$\otimes$\;
\ldots\;$\otimes$\;$\overline{x_k}$)\;=\;
$\overline{x_{\sigma(1)}}$\;$\otimes$\;\ldots
\;$\otimes$\;$\overline{x_{\sigma(k)}}$
\end{center}
for all $\overline{x_1}$\;$\in$\;$\mathcal{E}_1$,
$\overline{x_2}$\;$\in$\;$\mathcal{E}_2$, \ldots,
$\overline{x_k}$\;$\in$\;$\mathcal{E}_k$ is called  a $\sigma$-tensor permutation operator (TPO).\\
 If $n$\;=\;2, then say
that $U_{\sigma}$  is a tensor commutation operator.
\end{definition}
\begin{proposition}\label{tm32}
Consider the  $\mathbb{C}$- vector spaces $\mathcal{E}_1$,
$\mathcal{E}_2$, \ldots, $\mathcal{E}_k$, $\mathcal{F}_1$,
$\mathcal{F}_2$, \ldots, $\mathcal{F}_k$,   a permutation $\sigma$ of
\{${1, 2,\ldots , k}$\} and a $\sigma$-TPO\\ 
$U_{\sigma}\in \mathcal{L}\left(\mathcal{F}_1\otimes\ldots\otimes \mathcal{F}_k,
\mathcal{F}_{\sigma(1)}\otimes\ldots\otimes\mathcal{F}_{\sigma(k)}\right)$.
Then, for all $\phi_1\in \mathcal{L}\left(\mathcal{E}_1,\mathcal{F}_1\right)$,
$\phi_2\in \mathcal{L}\left(\mathcal{E}_2,\mathcal{F}_2\right)$, \ldots,
$\phi_k\in \mathcal{L}\left(\mathcal{E}_k,\mathcal{F}_k\right)$
\begin{equation*}
U_{\sigma}\cdot\left(\phi_1 \otimes  \ldots \otimes \phi_k\right)
=\left(\phi_{\sigma(1)}\otimes\ldots\otimes
\phi_{\sigma(k)}\right)\cdot V_{\sigma}
\end{equation*}
where $V_{\sigma}\in \mathcal{L}(\mathcal{E}_1\otimes\mathcal{E}_2\otimes\ldots\otimes \mathcal{E}_k,
\mathcal{E}_{\sigma(1)}\otimes\mathcal{E}_{\sigma(2)}\otimes\ldots\otimes\mathcal{E}_{\sigma(k)})$ is a $\sigma$-TPO.  
\end{proposition}

\begin{proof}
$\phi_1\otimes\ldots\otimes\phi_k
\in\mathcal{L}\left(\mathcal{E}_1\otimes
\ldots\otimes
\mathcal{E}_k,\mathcal{F}_1\otimes\ldots\otimes \mathcal{F}_k,
\right)$, thus\\
$U_{\sigma}\cdot(\phi_1\otimes\phi_2\otimes\ldots\otimes\phi_k)
\in\mathcal{L}\left(\mathcal{E}_1\otimes\ldots\otimes \mathcal{E}_k,
\mathcal{F}_{\sigma(1)}\otimes
\ldots\otimes\mathcal{F}_{\sigma(k)}\right)$.\\
$\left(\phi_{\sigma(1)}\otimes\phi_{\sigma(2)}\otimes\ldots\otimes\phi_{\sigma(k)}\right)\cdot V_{\sigma}
\in\mathcal{L}\left(\mathcal{E}_1\otimes\ldots\otimes \mathcal{E}_k,
\mathcal{F}_{\sigma(1)}\otimes\ldots\otimes\mathcal{F}_{\sigma(k)}\right)$\\
If $\overline{x_1}\in\mathcal{E}_1$, 
$\overline{x_2}\in\mathcal{E}_2$, \ldots,
$\overline{x_k}\in\mathcal{E}_k$,\\
$U_{\sigma}\cdot\left(\phi_1\otimes\ldots\otimes\phi_k\right)
\left(\overline{x_1}\otimes\ldots\otimes\overline{x_k}\right)\\
=U_{\sigma}\left[\phi_1\left(\overline{x_1}\right)\otimes\phi_2\left(\overline{x_2}\right)\otimes\ldots\otimes\phi_k\left(\overline{x_k}\right)\right]\\
=\phi_{\sigma(1)}\left(\overline{x_{\sigma(1)}}\right)\otimes\ldots\otimes\phi_{\sigma(k)}\left(\overline{x_{\sigma(k)}}\right)$\\
( Since $U_{\sigma}$ is a  TPO)\\
$=\left(\phi_{\sigma(1)}\otimes\ldots\otimes\phi_{\sigma(k)}\right)\left(\overline{x_{\sigma(1)}}\otimes\ldots
\otimes\overline{x_{\sigma(k)}}\right)\\=\left(\phi_{\sigma(1)}\otimes\ldots\otimes\phi_{\sigma(k)}\right)\cdot V_{\sigma}\left(\overline{x_1}\otimes\ldots\otimes\overline{x_k}\right)$.
\end{proof}
\begin{proposition}\label{prop6}
If  $U_{\sigma}$ is a $\sigma$-TPO, then
its transpose ${U_{\sigma}}^t$  is the ${\sigma}^{-1}$-TPO $U_{\sigma^{-1}}$.(Cf. for example $\cite{rak05}$)\\
\end{proposition}

\section*{Tensor permutation matrices}
\begin{definition} Consider the $\mathbb{C}$-vector spaces
$\mathcal{E}_1$, $\mathcal{E}_2$, \ldots, $\mathcal{E}_k$ of
dimensions $n_1$, $n_2$, \ldots, $n_k$ and $\sigma$-TPO
$U_{\sigma}$\;$\in\;$ $\mathcal{L}$($\mathcal{E}_1$\;$\otimes$\;
$\mathcal{E}_2$\;$\otimes$\;\ldots\;$\otimes$ $\mathcal{E}_k$,
$\mathcal{E}_{\sigma(1)}$\;$\otimes$\;
$\mathcal{E}_{\sigma(2)}$\;$\otimes$\;\ldots\;$\otimes$\;
$\mathcal{E}_{\sigma(k)}$).
Let\\
$\mathcal{B}_1=\left(\overline{e_{11}},\overline{e_{12}},\ldots,\overline{e_{1n_1}}\right)$ be a basis of
$\mathcal{E}_1$;\\
$\mathcal{B}_2=\left(\overline{e_{21}},\overline{e_{22}},\ldots,\overline{e_{2n_2}}\right)$ be a basis of
$\mathcal{E}_2$;\\\ldots\\
$\mathcal{B}_k=\left(\overline{e_{k1}},\overline{e_{k2}},\ldots,\overline{e_{kn_k}}\right)$ be a basis of
$\mathcal{E}_k$.\\
$\textbf{U}_{\sigma}$ the matrix of $U_{\sigma}$ with respect to the couple
of bases\\
\noindent $\left(\mathcal{B}_1\otimes\mathcal{B}_2\otimes \ldots \otimes \mathcal{B}_k, \mathcal{B}_{\sigma(1)}
\otimes\mathcal{B}_{\sigma(2)}\otimes \ldots\otimes\mathcal{B}_{\sigma(k)}\right)$. The square matrix $\textbf{U}_{\sigma}$ of
dimensions $n_1\times n_2\times\ldots \times n_k$ is independent of the
bases $\mathcal{B}_1$, $\mathcal{B}_2$,\ldots , $\mathcal{B}_k$.
Call this matrix a $\sigma$-TPM $n_1\otimes
n_2\otimes\ldots \otimes n_k$.
\end{definition}
According to the proposition \ref{prop6}, we have the following proposition.
\begin{proposition} 
A $\sigma$-TPM $\textbf{U}_{\sigma}$ is an orthogonal matrix, that is $\textbf{U}_{\sigma}^{-1}=\textbf{U}_{\sigma}^T$. 
\end{proposition} 
\begin{proposition}\label{prop7}
Let $\textbf{U}_{\sigma}$ be $\sigma$-TPM $n_1\otimes
n_2\otimes\ldots \otimes n_k$ and $\textbf{V}_{\sigma}$ a $\sigma$-TPM $m_1\otimes
m_2\otimes\ldots \otimes m_k$. Then, for all matrices
$\textbf{A}_1$, $\textbf{A}_2$,\ldots, $\textbf{A}_k$, of dimensions, respectively, $m_1\times n_1$,
$m_2\times n_2$, \ldots, $m_k\times n_k$\\
$\textbf{U}_{\sigma}\cdot\left(\textbf{A}_1\otimes 
\ldots\otimes\textbf{A}_k\right)\cdot\textbf{V}_{\sigma}^T=\textbf{A}_{\sigma(1)}\otimes
 \ldots\otimes
\textbf{A}_{\sigma(k)}$
\end{proposition}

\begin{proof} Let $A_1\in\mathcal{L}\left(\mathcal{E}_1,\mathcal{F}_1\right)$,
$A_2\in\mathcal{L}\left(\mathcal{E}_2,\mathcal{F}_2\right)$, \ldots,
$A_k\in\mathcal{L}\left(\mathcal{E}_k,\mathcal{F}_k\right)$. Their matrices with
respect to couple of bases $\left(\mathcal{B}_1,\mathcal{B}^{'}_1\right)$, $\left(\mathcal{B}_2,\mathcal{B}^{'}_2\right)$,\ldots , $\left(\mathcal{B}_k,\mathcal{B}^{'}_k\right)$
are respectively $\textbf{A}_1$, $\textbf{A}_2$,\ldots, $\textbf{A}_k$. Then,
$\textbf{A}_1\otimes \textbf{A}_2\otimes \ldots\otimes\textbf{A}_k$ is the matrix of
$A_1\otimes A_2\otimes\ldots\;\otimes A_k$ with respect to
$\left(\mathcal{B}_1\otimes\mathcal{B}_2\otimes
\ldots\otimes\mathcal{B}_k,\mathcal{B}{'}_1\otimes\mathcal{B}{'}_2\otimes
\ldots\otimes\mathcal{B}{'}_k\right)$.\\
But,  $A_1\otimes A_2\otimes \ldots\otimes A_k
\in\mathcal{L}\left(\mathcal{E}_1\otimes
\mathcal{E}_2\otimes\ldots\otimes
\mathcal{E}_k,\mathcal{F}_1\otimes\mathcal{F}_2\otimes\ldots\otimes \mathcal{F}_k
\right)$ and
$A_{\sigma(1)}\otimes A_{\sigma(2)}\otimes \ldots\otimes
A_{\sigma(k)}
\in \mathcal{L}\left(\mathcal{E}_{\sigma(1)}\otimes\;\ldots\otimes
\mathcal{E}_{\sigma(k)},\mathcal{F}_{\sigma(1)}\otimes\ldots\otimes\mathcal{F}_{\sigma(k)}\right)$,\\
 thus\\
  $U_{\sigma}\cdot\left(A_1\otimes  \ldots\otimes A_k\right)$,
$\left(A_{\sigma(1)}\otimes  \ldots\otimes
A_{\sigma(k)}\right)\cdot V_{\sigma}$
$\in \mathcal{L}\left(\mathcal{E}_1\otimes\ldots\otimes\mathcal{E}_k,
\mathcal{F}_{\sigma(1)}\otimes\ldots\otimes\mathcal{F}_{\sigma(k)}\right)$.\\
$\textbf{A}_{\sigma(1)}\otimes \textbf{A}_{\sigma(2)}\otimes \ldots\otimes
\textbf{A}_{\sigma(k)}$ is the matrix of $A_{\sigma(1)}\otimes
A_{\sigma(2)}\otimes \ldots\otimes A_{\sigma(k)}$ with respect to
$\left(\mathcal{B}_{\sigma(1)}\otimes \ldots\otimes\mathcal{B}_{\sigma(k)},\mathcal{B}{'}_{\sigma(1)}\otimes \ldots\otimes\mathcal{B}{'}_{\sigma(k)}\right)$. Thus $(\textbf{A}_{\sigma(1)}\otimes
 \ldots\otimes
\textbf{A}_{\sigma(k)})\cdot\textbf{V}_{\sigma}$ is the one of
$(A_{\sigma(1)}\otimes \ldots\otimes
A_{\sigma(k)})\cdot V_{\sigma}$ with respect to
$\left(\mathcal{B}_1\otimes\ldots\otimes\mathcal{B}_k, \mathcal{B}{'}_{\sigma(1)}\otimes
\ldots \otimes \mathcal{B}{'}_{\sigma(k)}\right)$ .\\
 $\textbf{U}_{\sigma}\cdot\left(\textbf{A}_1\otimes \textbf{A}_2\otimes
\ldots\otimes\textbf{A}_k\right)$ is the matrix of $U_{\sigma}\cdot\left(A_1\otimes
A_2\otimes \ldots\otimes A_k\right)$ with respect to the same couple of bases.\\
According to the proposition \ref{tm32},
we have\\
$\textbf{U}_{\sigma}\cdot\left(\textbf{A}_1\otimes 
\ldots \otimes\textbf{A}_k\right)
=\left(\textbf{A}_{\sigma(1)}\otimes
\ldots\otimes
\textbf{A}_{\sigma(k)}\right)\cdot\textbf{V}_{\sigma}$
\end{proof}
\begin{proposition}\label{tm41}
The matrix $\textbf{U}_{\sigma}$ is a $\sigma$-TPM $n_1\otimes
n_2\otimes\ldots \otimes n_k$ if, and only if, for all
$\textbf{a}_1\in\mathbb{C}^{n_1\times 1}$,
$\textbf{a}_2\in\mathbb{C}^{n_2\times 1}$,\ldots,
$\textbf{a}_k\in\mathbb{C}^{n_k\times 1}$\\
$\textbf{U}_{\sigma}\cdot(\textbf{a}_1\otimes \ldots\otimes
\textbf{a}_k)$\;=\;$\textbf{a}_{\sigma(1)}\otimes
\ldots\otimes \textbf{a}_{\sigma(k)}$

\end{proposition}

\begin{proof} $"\Longrightarrow"$  It is evident from the proposition \ref{prop7}.\\
$"\Longleftarrow"$ Suppose that for all\\
$\textbf{a}_1\in\mathbb{C}^{n_1\times 1}$,
$\textbf{a}_2\in\mathbb{C}^{n_2\times 1}$,\ldots,
$\textbf{a}_k\in\mathbb{C}^{n_k\times 1}$

$\textbf{U}_{\sigma}\cdot(\textbf{a}_1\otimes \ldots\otimes
\textbf{a}_k)$\;=\;$\textbf{a}_{\sigma(1)}\otimes
\ldots\otimes \textbf{a}_{\sigma(k)}$

Let $\overline{a_1}$\;$\in$\;$\mathcal{E}_1$,
$\overline{a_2}$\;$\in$\;$\mathcal{E}_2$, \ldots,
$\overline{a_k}$\;$\in$\;$\mathcal{E}_k$ and $\mathcal{B}_1$,
$\mathcal{B}_2$,\ldots , $\mathcal{B}_k$ be some bases of
$\mathcal{E}_1$, $\mathcal{E}_2$, \ldots, $\mathcal{E}_k$ where the
components of $\overline{a_1}$, $\overline{a_2}$,\ldots,
$\overline{a_k}$ form the unicolumn matrices $\textbf{a}_1$,
$\textbf{a}_2$,\ldots, $\textbf{a}_k$. The $\sigma$-TPO 
$U_{\sigma}$\;$\in\;$ $\mathcal{L}$($\mathcal{E}_1$\;$\otimes$\;
$\mathcal{E}_2$\;$\otimes$\;\ldots\;$\otimes$ $\mathcal{E}_k$,
$\mathcal{E}_{\sigma(1)}$\;$\otimes$\;
$\mathcal{E}_{\sigma(2)}$\;$\otimes$\;\ldots\;$\otimes$\;
$\mathcal{E}_{\sigma(k)}$) whose matrix with respect to
($\mathcal{B}_1$$\otimes$$\mathcal{B}_2$$\otimes$ \ldots$\otimes$
$\mathcal{B}_k$, $\mathcal{B}_{\sigma(1)}$$\otimes$
$\mathcal{B}_{\sigma(2)}$$\otimes$ \ldots$\otimes$
$\mathcal{B}_{\sigma(k)}$ ) is $\textbf{U}_{\sigma}$. Thus\\
$U_{\sigma}$($\overline{a_1}$\;$\otimes$\;\ldots\;$\otimes$\;$\overline{a_k}$)\;=\;
$\overline{a_{\sigma(1)}}$\;$\otimes$\;\ldots
\;$\otimes$\;$\overline{a_{\sigma(k)}}$.
This is true for all $\overline{a_1}$\;$\in$\;$\mathcal{E}_1$,
$\overline{a_2}$\;$\in$\;$\mathcal{E}_2$, \ldots,
$\overline{a_k}$\;$\in$\;$\mathcal{E}_k$.\\
Since $U_{\sigma}$ is a $\sigma$-TPO,
$\textbf{U}_{\sigma}$ is a $\sigma$-TPM $n_1\otimes
n_2\otimes\ldots \otimes n_k$ .
\end{proof}
The application of this proposition to any two unicolumn matrices leads us to the following remark. 
\begin{remark}\label{rmk9}
Consider the function $L$ from the set of all matrices into the set of all unicolumn matrices. For a $n\times p$-matrix $\textbf{X}=\begin{bmatrix}
  x_{11} & x_{12} & \ldots & x_{1p} \\
  x_{21} & x_{22} & \ldots & x_{2p}\\
  \ldots & \ldots & \ldots & \ldots  \\
  x_{n1} &  x_{n2} & \ldots & x_{np} \\
   \end{bmatrix}$,
   $L\left(\textbf{X}\right)=\begin{bmatrix}x_{11}\\
   x_{12}\\
   \vdots\\
   x_{1p}\\
   x_{21}\\
   x_{22}\\
   \vdots\\
   x_{2p}\\
   \vdots\\
   x_{n1}\\
   x_{n2}\\
   \vdots\\
   x_{np}\\
   \end{bmatrix}$.
    The relation $\textbf{U}_{n\otimes p}\cdot L\left(\textbf{X}\right)=L\left(\textbf{X}^T\right)$ can be obtained easily.
   
\end{remark}
\begin{example}
 Consider the matrix equation $\textbf{A}\cdot\textbf{X}\cdot\textbf{B}=\textbf{C}$, with respect to the unknown $\textbf{X}\in \mathbb{C}^{n\times q}$, where  $\textbf{A}\in \mathbb{C}^{m\times n}$, $\textbf{B}\in \mathbb{C}^{p\times q}$ and $\textbf{C}\in \mathbb{C}^{m\times q}$ are known. This equation can be transformed to the system of linear equations, whose matrix equation is $\cite{ikr77}$
\begin{equation}\label{eq1}
\left(\textbf{A}\otimes \textbf{B}^T\right)\cdot L\left(\textbf{X}\right)=L\left(\textbf{C}\right)
\end{equation}
 or equivalently
 \begin{equation}\label{eq2}
  \left(\textbf{B}^T\otimes \textbf{A}\right)\cdot L\left(\textbf{X}^T\right)=L\left(\textbf{C}^T\right)
\end{equation}  
   The equation $(\ref{eq2})$ can be obtained by multiplying the equation $(\ref{eq1})$ by the $m\otimes q$ TCM $\textbf{U}_{m\otimes q}$ and in using the proposition \ref{prop7} and the remark $\ref{rmk9}$.\\
    Mutually, the equation $(\ref{eq1})$ can be obtained by multiplying the equation $(\ref{eq2})$ by the $q\otimes m$ TCM $\textbf{U}_{q\otimes m}$.
 \end{example}
 \section*{Conclusion}
 We have generalized a property of TPM's. Two TPM's permutate tensor product of rectangle matrices. The example show the utility of the property. It suffices to transform a matrix linear  equation to a matrix linear equation of the form $\textbf{A}\textbf{X}=\textbf{B}$. Another matrix linear equation of this form can be deduced by using a TCM and by applying the generalization. 
 
\subsection*{Acknowledgments}
The author would like to thank Andriamifidisoa Ramamonjy of the Department of Mathematics and Informatics of the University of Antananarivo for his help in preparing the manuscript.

\end{document}